\documentclass{amsart}

\usepackage{amsmath,amssymb,amsthm,a4wide}
\usepackage{pgfplots}
\usepackage{graphicx,tikz}

\newtheorem*{thm}{Theorem}

\newtheorem*{corollary}{Corollary}

\theoremstyle{definition}

\theoremstyle{remark}

\begin{document}

\title[]{On the Spectral Resolution of Products\\ of Laplacian Eigenfunctions}

\author[]{Stefan Steinerberger}
\address{Department of Mathematics, Yale University, New Haven, CT 06511, USA}
\email{stefan.steinerberger@yale.edu}

\begin{abstract} We study products of eigenfunctions of the Laplacian $-\Delta \phi_{\lambda} = \lambda \phi_{\lambda}$ on compact manifolds. If $\phi_{\mu}, \phi_{\lambda}$ are two
eigenfunctions and $\mu \leq \lambda$, then one would perhaps expect their product $\phi_{\mu}\phi_{\lambda}$ to be mostly a linear combination of eigenfunctions with eigenvalue close to $\lambda$. This can faily quite dramatically: on $\mathbb{T}^2$, we see that
$$ 2\sin{(n x)} \sin{((n+1) x)} = \cos{(x)}  - \cos{( (2n+1) x)} $$
has half of its $L^2-$mass at eigenvalue 1. Conversely, the product
$$ \sin{(n x)} \sin{(m y)} \qquad \mbox{lives at eigenvalue} \quad \max{\left\{m^2,n^2\right\}} \leq m^2 + n^2 \leq 2\max{\left\{m^2,n^2\right\}}$$
and the heuristic is valid. We show that the main reason is that in the first example 'the waves point in the same direction':  if the heuristic fails and multiplication carries $L^2-$mass to lower frequencies, then
$\phi_{\mu}$ and $\phi_{\lambda}$ are strongly correlated at scale $ \sim \lambda^{-1/2}$ (the shorter wavelength)
$$ \left\| \int_{M}{ p(t,x,y)( \phi_{\lambda}(y) - \phi_{\lambda}(x))( \phi_{\mu}(y) - \phi_{\mu}(x)) dy} \right\|_{L^2_x} \gtrsim \| \phi_{\mu}\phi_{\lambda}\|_{L^2},$$
where $p(t,x,y)$ is the classical heat kernel and $t \sim \lambda^{-1}$. This turns out to be a fairly fundamental principle and is even valid for the Hadamard product of eigenvectors of a Graph Laplacian.
\end{abstract}

\maketitle

\section{Introduction and main results}
\subsection{Introduction.} We start by considering  Laplacian eigenfunctions $ -\Delta \phi_{\lambda} = \lambda \phi_{\lambda}$ on smooth, compact Riemannian manifolds with or without boundary (and Dirichlet or Neumann conditions in the case of a boundary). Ultimately, the results are so general that they will easily transfer to the case of finite graphs equipped with a Graph Laplacian. We consider a basic question.
\begin{quote}
\textbf{Question.} What is the behavior of the spectral resolution of the product $\phi_{\mu}\phi_{\lambda}$? 
\end{quote} 

Except for the special case where additional structure is present (see for example Sarnak \cite{Sarnak},  Bernstein \& Reznikoff \cite{bern}
and Kr\"otz \& Stanton \cite{krotz}), there seem to be very few results in the literature. Already the $L^2-$size of the product is a highly nontrivial quantity: a seminal result due to Burq-G\'{e}rard-Tzetkov \cite{burq}
establishes
$$ \|\phi_{\mu}\phi_{\lambda} \|_{L^2} \lesssim \min( \sqrt{\lambda}, \sqrt{\mu})^{} \| \phi_{\lambda}\|_{L^2} \|\phi_{\mu}\|_{L^2}$$
on compact two-dimensional manifolds without boundary (the same authors extended the result to higher dimensions \cite{burq2}). Similar results were established on two- and three-dimensional
manifolds with boundary by Blair, Smith \& Sogge \cite{blair}. The problem of obtaining more informative bounds on $\left\langle \phi_{\mu}\phi_{\lambda}, \phi_{\nu}\right\rangle$ is clearly fairly fundamental by itself but also
 occurs somewhat frequently
in applications (see e.g. \cite{fib, mas, mas2} for problems arising in dimensionality reduction of high-dimensional data or \cite{bronstein} for pattern matching of deformable shapes). The problem is trivial on the torus $\mathbb{T}^d$ but already on $\mathbb{S}^2$ the arising integrals involve Clebsch-Gordon coefficients and are deeply connected to concepts in representation theory.

\subsection{A Napkin Computation} This section contains a very heuristic derivation of how one would hope to describe the phenomenon. We shall assume that $\mu \leq \lambda$ throughout the paper and start with a basic observation: the range of Theorems that one could prove seems quite limited because the spectrum of $\phi_{\mu}\phi_{\lambda}$ can behave in very different ways. This is already easily seen on the one-dimensional torus $\mathbb{T}$.  The simple identity
$$ 2\sin{(n x)} \sin{((n+1) x)} =  \cos{(x)}  - \cos{( (2n+1) x)} $$
shows that two high-frequency eigenfunctions can easily conspire to have their product
located to a substantial extent at low frequencies. Lifting the example to $\mathbb{T}^2$ produces an analogous counterexample that should be contasted with the product
$$ \sin{(n x)} \sin{(m y)} \qquad \mbox{living at eigenvalue} \quad \max{\left\{m^2,n^2\right\}} \leq m^2 + n^2 \leq 2\max{\left\{m^2,n^2\right\}}.$$
We will now claim describe how these two examples can be understood as the outcome of 'waves moving transversally' and 'waves moving in parallel'.

\begin{figure}[h!]
\begin{tikzpicture}[scale=0.7]
\draw [black, thick,  domain=1:4, samples=40, smooth]  plot ({\x }, {sin(5*deg(\x))-0.5});
\draw [black, thick,  domain=1:4, samples=40,smooth]  plot ({\x }, {sin(5*deg(\x))} );
\draw [black, thick,  domain=1:4, samples=40,smooth]  plot ({\x }, {sin(5*deg(\x))+0.5});
\node at (4.5, 0) { * };
\draw [black, thick,  domain=-1.5:1.5, samples=40, smooth]  plot ( {6.5 + sin(5*deg(\x))-0.5}, { \x });
\draw [black, thick,  domain=-1.5:1.5, samples=40, smooth]  plot ( {7 + sin(5*deg(\x))-0.5}, { \x });
\draw [black, thick,  domain=-1.5:1.5, samples=40, smooth]  plot ( {7.5 + sin(5*deg(\x))-0.5}, { \x });
\node at (8.5, 0) { = };
\draw [black, thick,  domain=9:12, samples=40, smooth]  plot ({\x }, {sin(5*deg(\x))-0.5- 0.5*(\x-11)});
\draw [black, thick,  domain=9:12, samples=40,smooth]  plot ({\x }, {sin(5*deg(\x)) - 0.5*(\x-11)} );
\draw [black, thick,  domain=9:12, samples=40,smooth]  plot ({\x }, {sin(5*deg(\x))+0.5 - 0.5*(\x-11)});
\node at (12.75, 0) { + };
\draw [black, thick,  domain=-1.5:1.5, samples=40, smooth]  plot ( {14.5 +0.3*(\x)+ sin(5*deg(\x))-0.5}, { \x });
\draw [black, thick,  domain=-1.5:1.5, samples=40, smooth]  plot ( {15  +0.3*(\x) + sin(5*deg(\x))-0.5}, { \x });
\draw [black, thick,  domain=-1.5:1.5, samples=40, smooth]  plot ( {15.5 +0.3*(\x) + sin(5*deg(\x))-0.5}, { \x });
\end{tikzpicture}
\caption{Wave calculus: the product of generic transversal waves decouples into phase-shifts of the original functions without massive change in frequency.}
\end{figure}
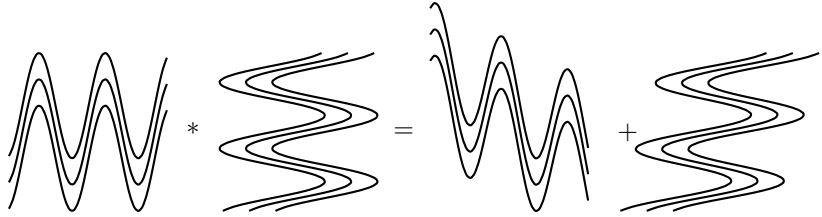
A very simple heuristic is that locally around a point $x_0$
$$ \phi_{\lambda}(x) \sim A \cos{\left(\left\langle x- x_0, \sqrt{\lambda}n_1 \right\rangle \right)} \qquad \mbox{and} \qquad  \phi_{\mu}(x) \sim B \cos{\left(\left\langle x- x_0, \sqrt{\mu}n_2 \right\rangle \right)},$$
where $n_1, n_2$ are two unit vectors. The product roughly behaves like
$$ \phi_{\lambda}(x) \phi_{\mu}(x) \sim \frac{AB}{2}\left(   \cos{\left(\left\langle x- x_0, \sqrt{\lambda}n_1 + \sqrt{\mu}n_2 \right\rangle \right)}  + \cos{\left(\left\langle x- x_0, \sqrt{\lambda}n_1 - \sqrt{\mu}n_2 \right\rangle \right)} \right).$$  
If $n_1, n_2$ are parallel (in the sense of having a large inner product), then some of the oscillation is moved to lower frequencies. Conversely, if they are transversal (having a small inner product) that neither
of the new terms is oscillating at a much lower frequency than $\lambda$. This suggests that the inner product $\left\langle n_1, n_2 \right\rangle $ should govern the amount with which energy transfer to lower frequencies is possible.
We now consider a sphere $\mathbb{S}_{x_0}$ centered at $x_0$ and having radius $\ll \lambda^{-1/2}$ and equipped with normalized Lebesgue measure $\sigma$. Then
$$ \int_{\mathbb{S}_{x_0}}{  ( \phi_{\lambda}(x) - \phi_{\lambda}(x_0)) ( \phi_{\mu}(x) - \phi_{\mu}(x_0)) d \sigma(x)} \sim AB \sqrt{\lambda \mu} \int_{\mathbb{S}_{x_0}}{ \left\langle n_1, x-x_0 \right\rangle \left\langle n_2, x-x_0 \right\rangle dx}.$$

\begin{figure}[h!]
\begin{tikzpicture}[scale=0.7]
\draw [black, thick,  domain=1:4, samples=40, smooth]  plot ({\x }, {sin(5*deg(\x))-0.5});
\draw [black, thick,  domain=1:4, samples=40,smooth]  plot ({\x }, {sin(5*deg(\x))} );
\draw [black, thick,  domain=1:4, samples=40,smooth]  plot ({\x }, {sin(5*deg(\x))+0.5});
\node at (4.5, 0) { * };
\draw [black, thick,  domain=5:8, samples=40, smooth]  plot ({\x }, {sin(5*deg(\x))-0.5});
\draw [black, thick,  domain=5:8, samples=40,smooth]  plot ({\x }, {sin(5*deg(\x))} );
\draw [black, thick,  domain=5:8, samples=40,smooth]  plot ({\x }, {sin(5*deg(\x))+0.5});
\node at (8.5, 0) { = };
\draw [black, thick,  domain=9:12, samples=40, smooth]  plot ({\x }, {sin(2*deg(\x))-0.5});
\draw [black, thick,  domain=9:12, samples=40,smooth]  plot ({\x }, {sin(2*deg(\x))} );
\draw [black, thick,  domain=9:12, samples=40,smooth]  plot ({\x }, {sin(2*deg(\x))+0.5});
\node at (12.75, 0) { + };
\draw [black, thick,  domain=13.5:16, samples=40, smooth]  plot ({\x }, {sin(8*deg(\x))-0.5});
\draw [black, thick,  domain=13.5:16, samples=40,smooth]  plot ({\x }, {sin(8*deg(\x))} );
\draw [black, thick,  domain=13.5:16, samples=40,smooth]  plot ({\x }, {sin(8*deg(\x))+0.5});
\end{tikzpicture}
\caption{Wave calculus: multiplying parallel waves leads to frequency change.}
\end{figure}
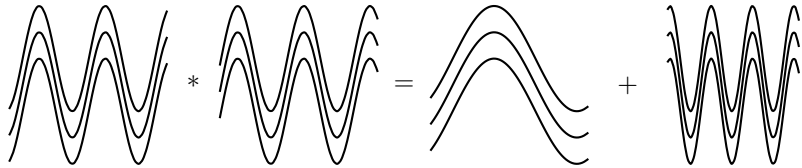

A simple computation shows that 
$$ \int_{\mathbb{S}_{x_0}}{ \left\langle n_1, x-x_0 \right\rangle \left\langle n_2, x - x_0\right\rangle dx} = \cos{\alpha}\int_{\mathbb{S}_{x_0}}{x_1^2 d\sigma(x)} = c_{\dim \mathbb{S}_{x_0}}  \cos{\alpha}.$$
This means that this local quantity will be small in absolute value if and only if the directions are close to transversal; in that case we expect little to no change for the location of the bulk of the spectrum.
Conversely, if the quantity is quite large, then it is conceivable for multiplication to allow for energy transfer to lower frequencies. We now average this operation over the entire manifold and
have as \textit{very} rough approximation 
$$  \int_{M}{ \left( \int_{\mathbb{S}_{x_0}}{  ( \phi_{\lambda}(x) - \phi_{\lambda}(x_0)) ( \phi_{\mu}(x) - \phi_{\mu}(x_0)) d \sigma(x)} \right)^2 d x_0}  \sim \int_{M}{  \phi_{\mu}^2 \phi_{\lambda}^2 \cos^2{ \angle \left(\frac{\nabla \phi_{\mu}}{|\nabla \phi_{\mu}|} \cdot   \frac{\nabla \phi_{\lambda}}{|\nabla \phi_{\lambda} |}\right)}  dx}.$$
The prediction from this heuristic is that if the left-hand side quantity is 'small' then the spectrum localizes essentially around frequency $\lambda$ whereas if the quantity is 'large', then there is substantial local
correlation of the gradient and mass transport to frequencies $\ll \lambda$ becomes possible.
The contribution of this short paper is to explain that this principle, 
 'waves pointing in the same direction' versus 'waves pointing different directions', can actually be made rigorous.

\subsection{Local correlation.}
To make this intuition precise, we will make use of the heat kernel $p(t,x,y):\mathbb{R}_{\geq 0} \times M \times M \rightarrow \mathbb{R}_{\geq 0}$. We note that for $(t,x) \in \mathbb{R}_{\geq 0} \times M$ fixed, we have
$$ p(t,x,y) \sim \begin{cases} t^{-d/2} \qquad &\mbox{if}~d(x,y) \lesssim \sqrt{t} \\ 0 \qquad &\mbox{otherwise} \end{cases} \qquad \mbox{and} \qquad \int_{M}{p(t,x,y) dy} = 1.$$
This suggests that one could think of $p(t,x,\cdot)$ as a local probability measure centered at $x$ and at scale $\sim t^{-1/2}$. Inspired by the heuristic in the previous section, a natural measure of local correlation in a point $x$ at scale $\sim t^{-1/2}$ is thus given by
$$ \mbox{local correlation} = \int_{M}{ p(t,x,y)( \phi_{\lambda}(y) - \phi_{\lambda}(x))( \phi_{\mu}(y) - \phi_{\mu}(x)) dy}.$$
Locally, it relates to $\cos{\angle (\nabla \phi_{\mu}, \nabla \phi_{\lambda})}$: we are interested in whether the quantity is close to 0 or not but it being $+1$ or $-1$ does not make a difference. The natural global
object is therefore given by its $L^2-$norm and we define our global notion of correlation as such.
$$ \mbox{Global correlation} =  \left\|  \int_{M}{ p(t,x,y)( \phi_{\lambda}(y) - \phi_{\lambda}(x))( \phi_{\mu}(y) - \phi_{\mu}(x)) dy} \right\|_{L^2_x}$$
should be large if and only if directions of oscillation tend to line up (leading to an energy transfer to low frequencies) and should only be small
if the directions of oscillation are mainly transversal (leading to no substantial change in the spectrum).
Before proceeding further, we would certainly
like to understand how the global correlation behaves for typical 'generic' eigenfunctions $\phi_{\lambda}$ and $\phi_{\mu}$.
A suitable candidate is given by arithmetic random waves on $\mathbb{T}^2$ that are believed to behave very much like typical 'chaotic' eigenfunctions \cite{berry,hejhal}. 
We thus fix two functions
$$ \phi_{\mu} = \sum_{m \in \mathbb{Z}^2 \atop |m|^2 = \mu}{ a_m e^{2\pi i \left\langle m, x \right\rangle}} \qquad \mbox{and} \qquad   \phi_{\lambda} = \sum_{n \in \mathbb{Z}^2 \atop |n|^2 = \lambda}{ b_n e^{2\pi i \left\langle n, x \right\rangle}}.$$
A somewhat lengthy but uninspired computation shows that $\mu \ll \lambda$ leads to
$$ \mbox{global correlation}  \sim  \left(1- e^{-t \mu} \right) \| \phi_{\mu}\phi_{\lambda}  \|_{L^2} $$
This immediately suggests a natural time-scale. If $\mu \ll \lambda$, then $t = \lambda^{-1}$ sets the scale for the local correlation functional at the shorter wavelength 
 and we would expect that the
$$ \mbox{global correlation at time } t=\frac{1}{\lambda} \sim  \left(1- e^{-\mu/\lambda} \right) \|  \phi_{\mu}\phi_{\lambda}   \|_{L^2} \sim \frac{\mu}{\lambda} \|  \phi_{\mu}\phi_{\lambda}  \|_{L^2} \qquad \mbox{is small.}$$
While this derivation used Fourier series, we emphasize that the prediction can be expected to be valued for generic chaotic eigenfunctions since global correlation is defined by local behavior.
Either the prediction carries over to the universal case \textit{or} the random wave model does not actually hold to the required extent; both cases would be of great interest but we are naturally inclined to believe
that the random wave model is accurate.

\subsection{Main result.} We can now state the main result. There are many possible ways of phrasing it: we state it as an algebraic fact. There exists
one unique point in time, slightly larger than $\lambda^{-1}$ such that the heat evolution of the product $\phi_{\mu}\phi_{\lambda}$ is \textit{exactly} given
by the local correlation function.

\begin{thm}[Main result] There exists a unique time $t \sim \log{(e\lambda/\mu)} \lambda^{-1}$ such that
 $$ \forall~x\in M \qquad   \left[e^{t \Delta} (\phi_{\mu}\phi_{\lambda})\right](x)  = \int_{\Omega}{ p(t,x,y)( \phi_{\lambda}(y) - \phi_{\lambda}(x))( \phi_{\mu}(y) - \phi_{\mu}(x)) dy}.$$
\end{thm}
Moreover, the unique time $t$ is the solution of the equation
$$ e^{-t \lambda} + e^{-t \mu} = 1.$$
The main implication for the spectral resolution of the product $\phi_{\mu}\phi_{\lambda}$ is as follows.
\begin{enumerate}
\item If the product $\phi_{\mu}\phi_{\lambda}$ has a nontrivial amount of $L^2-$mass at low frequencies $\lesssim \lambda/\log{(e\lambda/\mu)}$, then the heat equation has very little effect on these components
up to time $t \sim \log{(e\lambda/\mu)} \lambda^{-1}$ and we would expect $\| e^{t \Delta} \phi_{\mu} \phi_{\lambda} \|_{L^2}$ not to decay by much.
\item The main result implies that at time $t \sim \log{(e\lambda/\mu)} \lambda^{-1}$ the solution of the heat equation coincides with the local correlation functional. This implies that the global correlation functional is large  and thus there are large local correlations at small scale.
\end{enumerate}

We now explain this approach more formally. If $\mu \ll \lambda$ and the global interaction functional is small, then we expect the spectrum of the product to be mostly localized around frequency
$\sim \lambda$. Substantial energy transfer to lower frequencies does not occur.

\begin{corollary}[Transversal waves] If $\phi_{\lambda}, \phi_{\mu}$ satisfy, for $t$ solving $e^{-t \lambda} + e^{-t \mu} = 1$,
$$  \left\|  \int_{M}{ p(t,x,y)( \phi_{\lambda}(y) - \phi_{\lambda}(x))( \phi_{\mu}(y) - \phi_{\mu}(x)) dy} \right\|_{L^2_x}  \leq \delta \|\phi_{\mu}\phi_{\lambda} \|_{L^2},$$
then the amount of energy being moved to low frequencies is small 
$$      \sum_{\lambda_k \leq \lambda/\log{(e\lambda/\mu)}}^{}{  |\left\langle \phi_{\mu}\phi_{\lambda}, \phi_{k} \right\rangle|^2} \leq \delta^2 e^2 \| \phi_{\mu}\phi_{\lambda} \|^2_{L^2}.$$
\end{corollary}

If $\phi_{\mu}$ and $\phi_{\lambda}$ behave like random waves, then we expect the assumption to hold for $\delta \sim \mu / \lambda$.
There is also an inverse result: if there is a massive amount of correlation, then the energy cannot move to substantially higher frequencies: some of it
has to remain close to $\lambda$.

\begin{corollary}[Parallel waves] If $\phi_{\lambda}, \phi_{\mu}$ satisfy, for $t$ solving $e^{-t \lambda} + e^{-t \mu} = 1$,
$$  \left\|  \int_{M}{ p(t,x,y)( \phi_{\lambda}(y) - \phi_{\lambda}(x))( \phi_{\mu}(y) - \phi_{\mu}(x)) dy} \right\|_{L^2_x}  \geq \delta \| \phi_{\mu} \phi_{\lambda} \|_{L^2},$$
then there is at least some energy at frequencies $\lesssim \lambda \log{(1/\delta)}$
$$ \sum_{\lambda_k \leq  \lambda \log{(1/\delta)} }^{\infty}{|\left\langle \phi_{\mu}\phi_{\lambda}, \phi_{k} \right\rangle|^2}  \lesssim \delta^2  \| \phi_{\mu} \phi_{\lambda} \|^2_{L^2}.$$
\end{corollary}
We believe that these results raise many natural questions; is the spectrum of $\phi_{\mu} \phi_{\lambda}$ for $\mu \leq \lambda$
in the quantum chaotic regime mostly supported on frequencies in $(\lambda/\log{\lambda}, \lambda \log{\lambda})$?

\subsection{Graphs and Graph Laplacians.} The phenomenon is so fundamental that a version of it exists on finite graphs. Let $G=(V,E)$ be a finite, undirected, unweighted, connected graph on $|V| = n$
vertices. Functions are given by maps $f:V \rightarrow \mathbb{R}$ and are represented by elements in $\mathbb{R}^n$. We construct a Laplacian $\mathcal{L}: \mathbb{R}^n \rightarrow \mathbb{R}^n$ acting on functions via
$$   \left[\mathcal{L}f\right](v) = \sum_{w \sim_{E} v}{(f(w) - f(v))}.$$
This linear map corresponds to a symmetric positive-semidefinite matrix whose eigenvalues we denote by $0 = \lambda_1 < \lambda_2 \leq \lambda_3 \leq \dots \leq \lambda_n$ and eigenvectors by $\phi_1, \dots, \phi_n$. The definition of a Laplacian
enables us to define the heat flow $e^{t\Delta}$ via
$$  e^{t\Delta} f = \sum_{k=1}^{n}{ e^{-\lambda_k t } \left\langle f, \phi_k\right\rangle \phi_k} \quad \mbox{satisfying} \quad \partial_t e^{t\Delta} f  = \mathcal{L} e^{t\Delta} f $$
as desired. Being able to solve a heat equation allows the construction of a heat kernel $p(t,u,v) : \mathbb{R}_{\geq 0} \times \left\{1, \dots, n\right\} \times \left\{1, \dots, n\right\} \rightarrow \mathbb{R}_{\geq 0}.$
We use $\delta_{v}$ for the indicator function of the vertex $v$ and set
$$ p(t,u,v) = \left[e^{t\Delta} \delta_u\right](v).$$
It is easily seen that $p(t,u,v) = p(t,v,u)$. 
Moreover, by linearity of the heat equation,
$$ \left[e^{t\Delta}f\right](u) = \sum_{k=1}^{n}{p(t,u,v)f(v)}.$$
Alternatively, one could define the matrix evolution $e^{t \mathcal{L}}$ as a semigroup acting on matrices without appealing to the spectrum (we refer to Chung \cite{chung} for an introduction into
spectral graph theory).
We can now phrase our main result in the framework of Graphs.

\begin{thm}[Main result on Graphs] If $t$ is the unique solution of $e^{-\lambda t} + e^{-\mu t}=1$, then
 $$   e^{t \Delta} (\phi_{\mu}\phi_{\lambda})(u) = \sum_{v}^{}{ p(t,u,v)( \phi_{\lambda}(v) - \phi_{\lambda}(u))( \phi_{\mu}(v) - \phi_{\mu}(u)) dy}.$$
\end{thm}

We emphasize that the validity on Graphs show it to be a fundamental spectral principle; moreover, they are suitable toy examples because
all the computations reduce to basic linear algebra.

\section{Numerical Examples}
 We consider three examples (Fig. 3-5):  the Faulkner-Younger graph on 44 vertices, the Thomassen graph on 94 vertices and random graphs. 
These examples confirm the main principle: a large correlation at a local scale can cause a shift to lower frequencies whereas little correlation implies that the oscillation move transversally
and have their bulk continuing to oscillate at frequency $\sim \lambda$. While the notion of 'waves moving transversally' on a Graph with 44 vertices might seem like a bit of a stretch, 
the principle remains valid.

\subsection{Example 1: Faulkner-Younger Graph}
The Faulkner-Younger Graph on 44 vertices is a counterexample \cite{faulkner} (though not the first one) to an old 1880 conjecture of Tait \cite{tait}.
We compute the global correlation between all pairs of eigenvectors for $t$ solving $e^{-t\lambda} + e^{-t \mu} = 1$ and normalize by dividing by $\|\phi_{\mu}\phi_{\lambda}\|_{\ell^2}$. Global correlations are clustered around
$\sim 0.5$.
Notable outliers are
\begin{figure}[h!]
\includegraphics[width=0.3\textwidth]{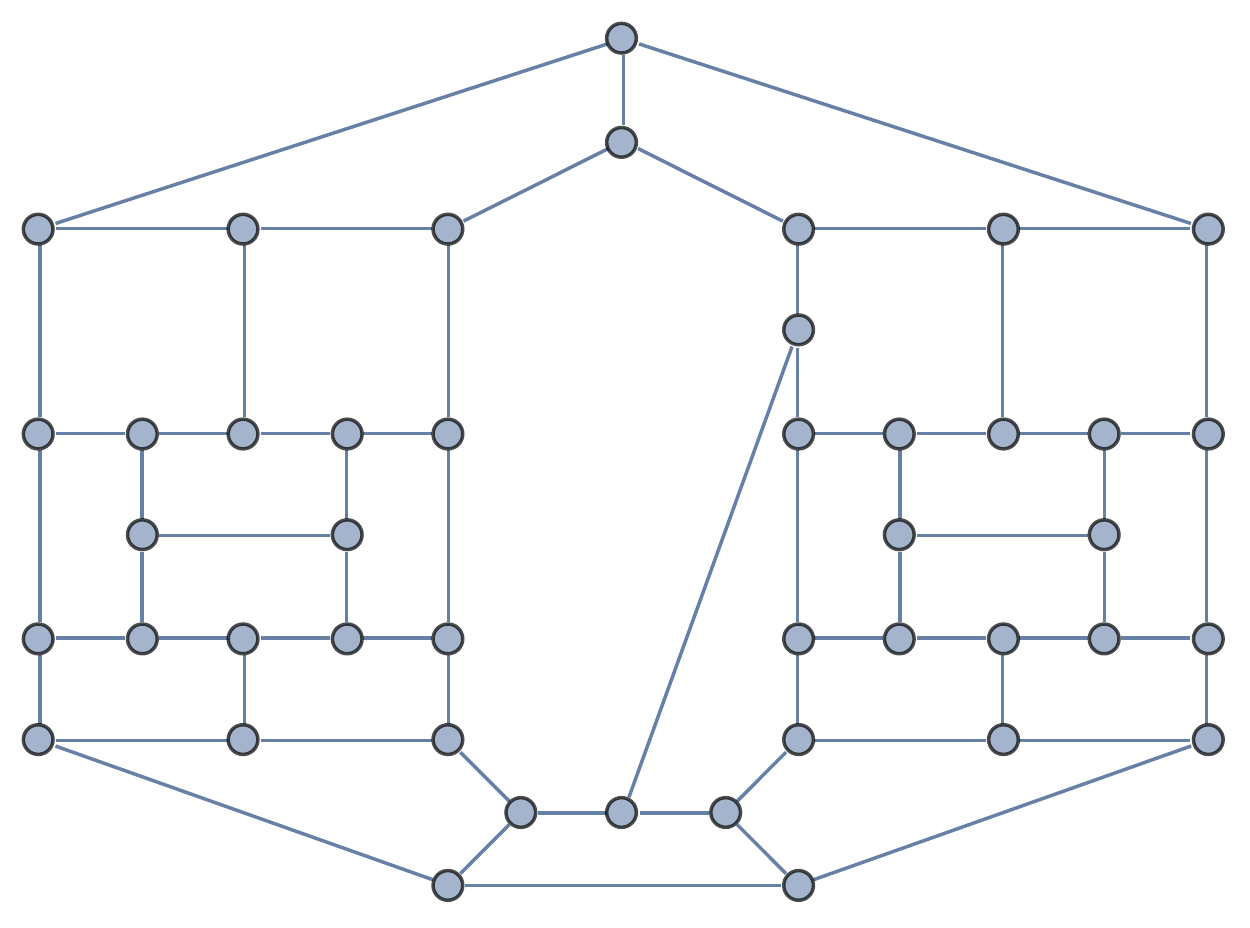}
\caption{The Faulkner-Younger Graph on 44 vertices.}
\end{figure}
\begin{enumerate}
\item $\phi_{12}, \phi_{15}$ ( $\sim 0.32$):  $ \phi_{12} \cdot \phi_{15}$ has $20\%$ of $\ell^2-$mass on $\phi_{15}$ and $58\%$ of $\ell^2-$mass on $\phi_{32}$.
\item $\phi_{12}, \phi_{24}$ ( $\sim 0.38$): $\phi_{12} \cdot \phi_{24}$ has $24\%$ of $\ell^2-$mass on $\phi_{26}$, $44\%$ on the highest 6 frequencies.
\item $\phi_{18}, \phi_{25}$ ( $\sim 0.7$): $\phi_{18} \cdot \phi_{25}$ has 60\% of $\ell^2-$mass on the 10 lowest frequencies.
\item $\phi_{40}, \phi_{43}$ ( $\sim 0.75$): $\phi_{40} \cdot \phi_{43}$ has 93\% of $\ell^2-$mass on the 25 lowest frequencies.
\item $\phi_{42}, \phi_{44}$ ( $\sim 0.74$): $\phi_{42} \cdot \phi_{44}$ has 94\% of $\ell^2-$mass on the 25 lowest frequencies.
\end{enumerate}

\subsection{Example 2: Thomassen Graph.}
The Thomassen Graph is the smallest cubic planar hypohamiltonian graph in the infinite family constructed by Thomassen \cite{th} to answer
a question of Chv\'{a}tal \cite{ch}; we expect eigenvectors to be structured but not arithmetically related. 

\begin{figure}[h!]
\includegraphics[width=0.5\textwidth]{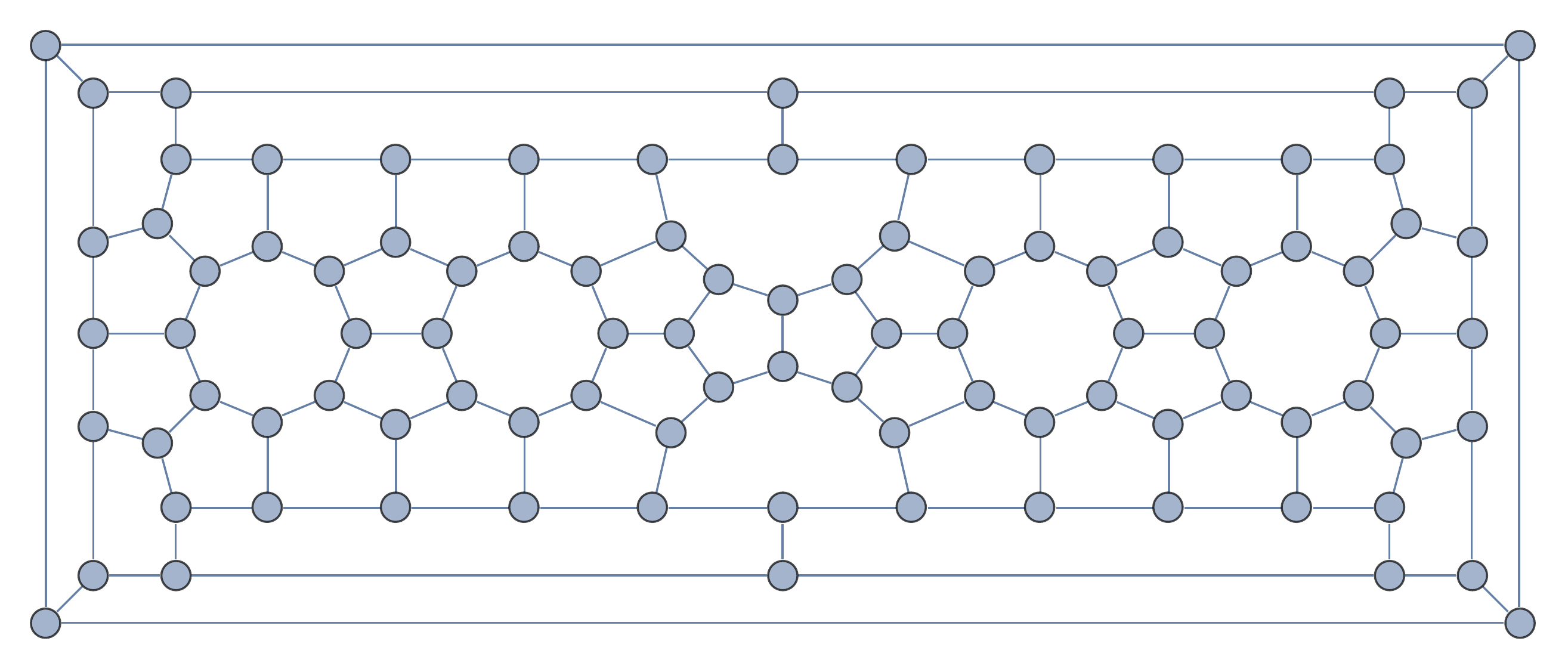}
\caption{The Thomassen Graph on 94 vertices.}
\end{figure}

We again compute the global correlation between all eigenvectors $\phi_{\mu}, \phi_{\lambda}$ and normalize by dividing by $\|\phi_{\mu}\phi_{\lambda}\|_{\ell^2}$. Most values at centered around $\sim 0.5$, notable outliers are

\begin{enumerate}
\item $\phi_{10}, \phi_{15}$ ( $\sim 0.14$): $\phi_{10} \cdot \phi_{15}$ has 78\% of $\ell^2-$mass on $\phi_{39}$.
\item $\phi_{23}, \phi_{43}$ ( $\sim 0.27$): $ \phi_{23} \cdot \phi_{43}$ has $15\%$ of mass on $\phi_{43}$, $14\%$ on $\phi_{90}$ and $56\%$ on $\phi_{91}$.
\item $\phi_{90}, \phi_{92}$ ( $\sim 0.75$): $\phi_{90} \cdot \phi_{92}$ has 76\% of $\ell^2-$mass on the 45 lowest frequencies.
\item $\phi_{90}, \phi_{93}$ ( $\sim 0.73$): $\phi_{90} \cdot \phi_{93}$ has 62\% of $\ell^2-$mass on $\phi_{40}$.
\item $\phi_{92}, \phi_{93}$ ( $\sim 0.93$): $\phi_{92} \cdot \phi_{93}$ has 57\% on $\phi_{4}$, 26\% on $\phi_{12}$ and 12\% on $\phi_{22}$.
\end{enumerate}

These results also indicate an unrelated question: is it common for Graph Laplacians to have two high-frequency eigenvectors whose Hadamard product lies at substantially lower frequencies? A possible
explanation would be that very high frequencies can only be obtained by essentially oscillating across every edge: multiplying two such objects would cancel the oscillation. This
fails for random/expander graphs so presumably some restriction on the geometry is required.

\subsection{Example 3:  Random Graphs.}
We quickly compare these results to are random graphs. It is easily observed that all the effects weaken and things become more uniform as the number of edges increase: this is not surprising since adding edges turns the graph into a more uniform object whose spectrum is more and more strongly localized.

\begin{figure}[h!]
\begin{minipage}[l]{.49\textwidth}
\begin{center}
\includegraphics[width = 0.8\textwidth]{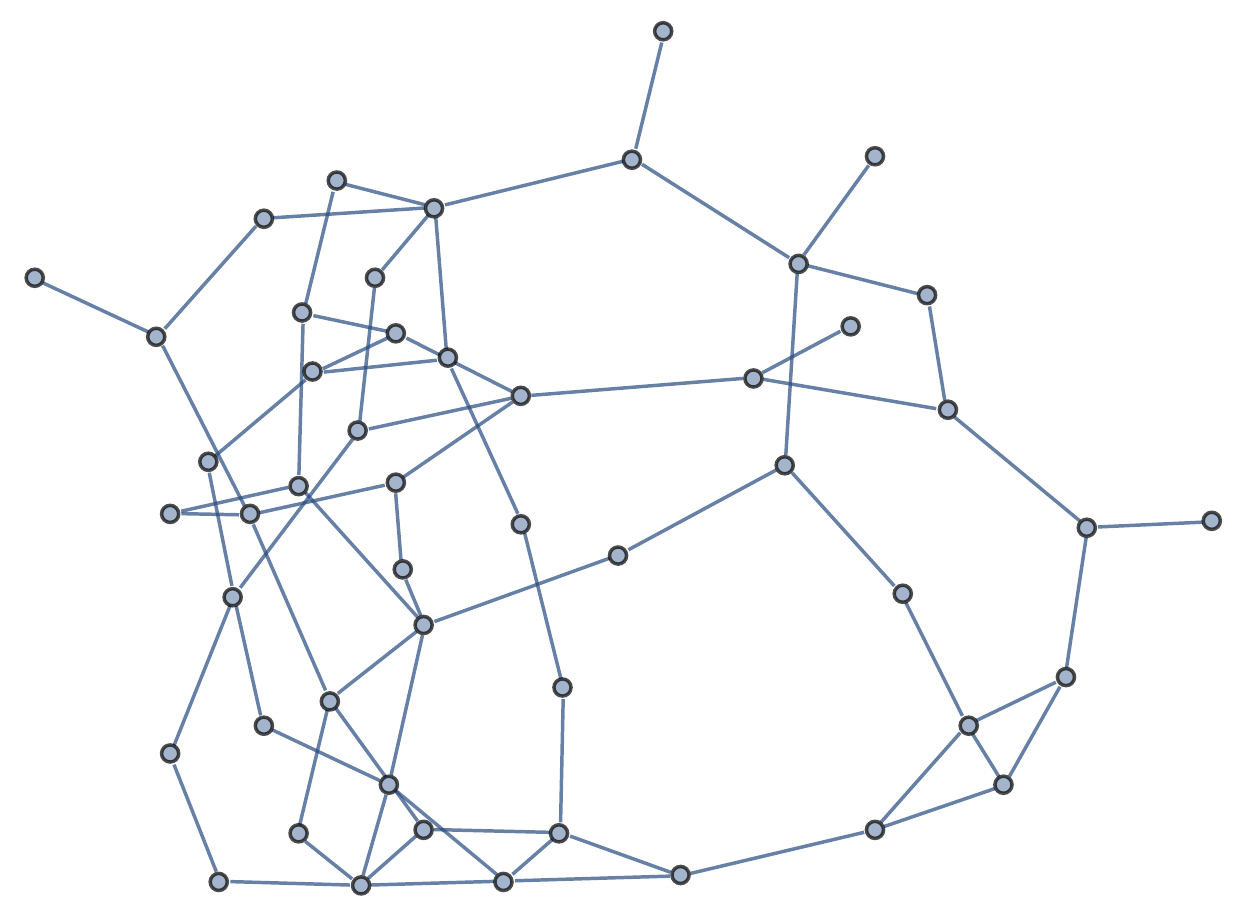}
\end{center}
\end{minipage} 
\begin{minipage}[r]{.49\textwidth}
\begin{center}
\includegraphics[width = 0.6\textwidth]{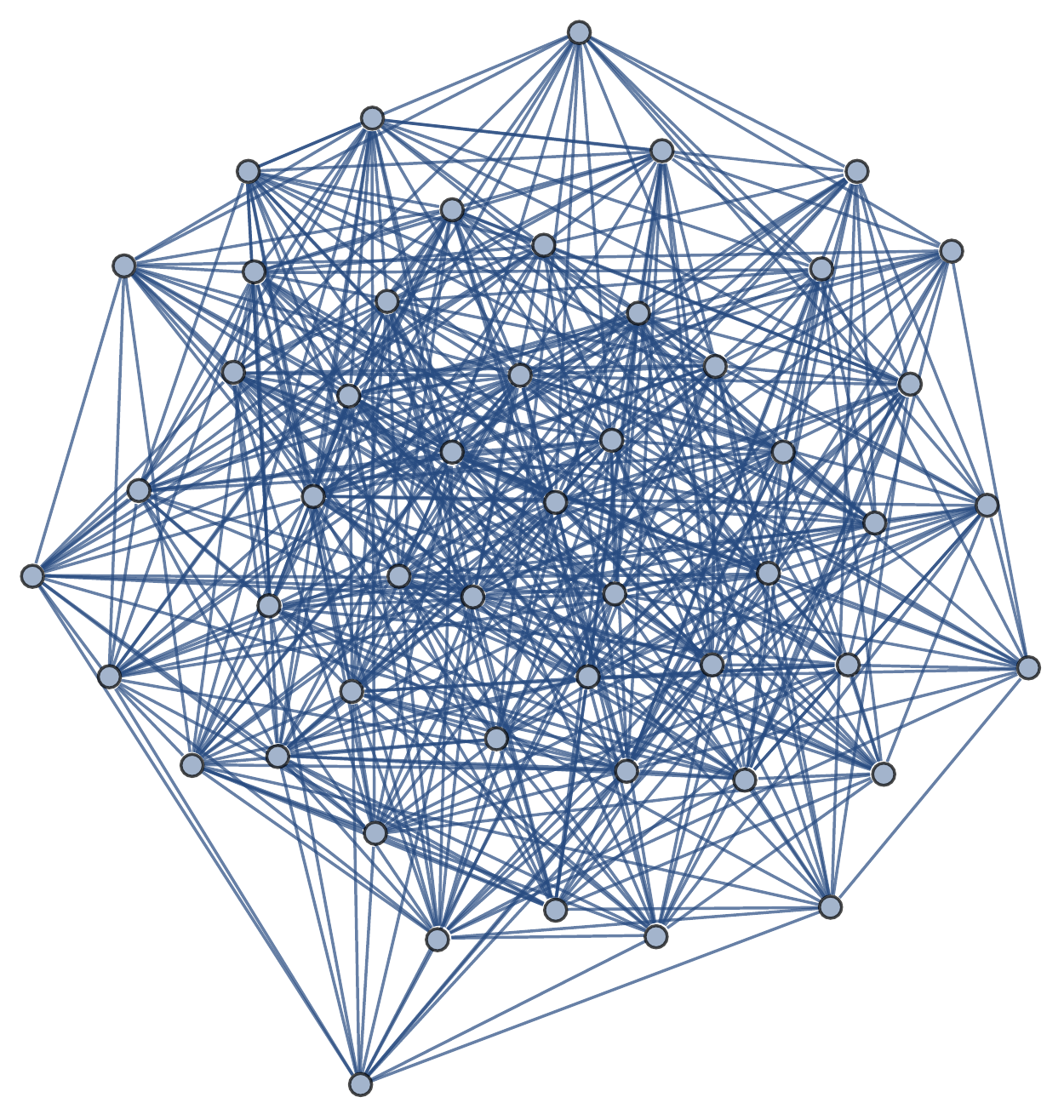}
\end{center}
\end{minipage} 
\caption{Random Graphs on 50 vertices and 70 edges (left) as well as 500 edges (right).}
\end{figure}

The typical correlation between any pair of eigenvectors becomes more uniform: in the example in Figure 4 (left), $\phi_{14} \cdot \phi_{15}$ ($\sim 0.64$) has 47\% of $\ell^2-$mass on the lowest 10 frequencies. $\phi_{15} \cdot \phi_{31}$  $\sim 0.73$) has 61\% on the first 9 frequencies. Moreover, $\phi_{48} \cdot \phi_{49}$
has 41\% on $\phi_{48}, \phi_{49}$. The correlation between eigenvectors of two randomly generated graphs (using Mathematica's RandomGraph) on 50 vertices and 100 edges has mean 0.48 and standard
deviation 0.11. In comparison, for 50 vertices and 1000 edges, the mean correlation is 0.50 with standard deviation 0.01.

\section{Proof}
\subsection{Proof of Theorem 1}
\begin{proof}  We phrase the proof in the language of diffusion since it dramatically simplifies notation. $\omega_x(t)$ refers to Brownian motion started in $x$ at time $t$. The main identity (see e.g. \cite{stein})
is the Feynman-Kac formula 
\begin{align*}
\left[e^{t \Delta} \phi_{\lambda}\right](x) = \mathbb{E}_{\omega}(\phi_{\lambda}(\omega_x(t))).
\end{align*}
Moreover, since Laplacian eigenfunctions diagonalize the heat flow
$$ \mathbb{E}_{\omega}(\phi_{\lambda}(\omega_x(t))) = e^{-\lambda t} \phi_{\lambda}(x).$$
We now start with a series of algebraic manipulations; adding 0 yields
\begin{align*}
\left[e^{t \Delta} \phi_{\mu}\phi_{\lambda}\right](x) - \phi_{\lambda}(x)\phi_{\mu}(x) &= \left[e^{t \Delta} \phi_{\mu}\phi_{\lambda}\right](x) - \phi_{\lambda}(x)\left[e^{t \Delta} \phi_{\mu}\right](x)\\
&+ \phi_{\lambda}(x)\left[e^{t \Delta} \phi_{\mu}\right](x) - \phi_{\lambda}(x)\phi_{\mu}(x).
\end{align*}
This can be rewritten as
\begin{align*}
 \left[e^{t \Delta} \phi_{\mu}\phi_{\lambda}\right](x) - \phi_{\lambda}(x)\phi_{\mu}(x) &= \mathbb{E}_{x}\left( \phi_{\lambda}(\omega_x(t))\phi_{\mu}(\omega_x(t))   \right)  -  \mathbb{E}_{x}\left( \phi_{\lambda}(x)\phi_{\mu}(\omega_x(t))   \right) \\
&+   \mathbb{E}_{x}\left( \phi_{\lambda}(x)\phi_{\mu}(\omega_x(t))   \right) -  \phi_{\lambda}(x) \phi_{\mu}(x)
\end{align*}
Using linearity of expectation
$$ \mathbb{E}_{x}\left( \phi_{\lambda}(\omega_x(t))\phi_{\mu}(\omega_x(t))   \right)  -  \mathbb{E}_{x}\left( \phi_{\lambda}(x)\phi_{\mu}(\omega_x(t))   \right) =
\mathbb{E}_{x}\left( \phi_{\mu}(\omega_x(t)) \left(\phi_{\lambda}(\omega_x(t))   - \phi_{\lambda}(x) \right)\right)$$
and
$$  \mathbb{E}_{x}\left( \phi_{\lambda}(x)\phi_{\mu}(\omega_x(t))   \right) -  \phi_{\lambda}(x) \phi_{\mu}(x) = \phi_{\lambda}(x) \mathbb{E}_{x}\left( \phi_{\mu}(\omega_x(t)) - \phi_{\mu}(x)  \right).$$ 
Writing $\phi_{\mu}(\omega_x(t)) = \phi_{\mu}(\omega_x(t)) - \phi_{\mu}(x)+ \phi_{\mu}(x)$ and plugging into the first term yields
\begin{align*}
 \left[e^{t \Delta} \phi_{\mu}\phi_{\lambda}\right](x) - \phi_{\lambda}(x)\phi_{\mu}(x) &= 
\phi_{\lambda}(x) \mathbb{E}_{x}\left( \phi_{\mu}(\omega_x(t)) - \phi_{\mu}(x)  \right)\\
& + \phi_{\mu}(x) \mathbb{E}_{x}\left( \phi_{\lambda}(\omega_x(t)) - \phi_{\lambda}(x)  \right) \\
&+ \mathbb{E}_{x}\left(   ( \phi_{\mu}(\omega_x(t)) - \phi_{\mu}(x)) ( \phi_{\lambda}(\omega_x(t)) - \phi_{\lambda}(x))      \right).
\end{align*}
However, two out of three of these expectations can be explicitly computed and this yields
\begin{align*}
 \left[e^{t \Delta} \phi_{\mu}\phi_{\lambda}\right](x) - \phi_{\lambda}(x)\phi_{\mu}(x) &= 
\phi_{\lambda}(x)( e^{-\mu t} - 1) \phi_{\mu}(x)  \\
& + \phi_{\mu}(x) ( e^{-\lambda t} - 1)  \phi_{\lambda}(x)   \\
&+ \mathbb{E}_{x}\left(   ( \phi_{\mu}(\omega_x(t)) - \phi_{\mu}(x)) ( \phi_{\lambda}(\omega_x(t)) - \phi_{\lambda}(x))      \right).
\end{align*}
At the same time, the expectation can be rewritten in terms of the heat kernel as
$$  \mathbb{E}_{x}\left(   ( \phi_{\mu}(\omega_x(t)) - \phi_{\mu}(x)) ( \phi_{\lambda}(\omega_x(t)) - \phi_{\lambda}(x))      \right) = \int_{\Omega}{ p(t,x,y)( \phi_{\lambda}(y) - \phi_{\lambda}(x))( \phi_{\mu}(y) - \phi_{\mu}(x)) dy}.$$
Summarizing, we obtain the fundamental identity
\begin{align*}
  \left[e^{t \Delta} \phi_{\mu}\phi_{\lambda}\right](x)  &= \left( e^{-\lambda t} + e^{-\mu t} - 1\right) \phi_{\lambda}(x) \phi_{\mu}(x) \\
&+  \int_{\Omega}{ p(t,x,y)( \phi_{\lambda}(y) - \phi_{\lambda}(x))( \phi_{\mu}(y) - \phi_{\mu}(x)) dy}.
\end{align*}
It remains to determine the size of $t$ such that
$$e^{-\lambda t} + e^{-\mu t} = 1.$$
It is clear that the left-hand side is bigger than 1 for $t$ sufficiently small. 
Setting $t = \alpha \log{(e\lambda/\mu)} \lambda^{-1}$ shows that
$$ e^{-\lambda t} + e^{-\mu t}  = \left(  \frac{e \mu}{\lambda} \right)^{\alpha} +  \left(  \frac{e \mu}{\lambda} \right)^{\alpha\frac{\mu}{\lambda} }.$$
A computation shows that $(ex)^{0.8} + (ex)^{0.8x} \geq 1$ on $(0, e)$. This shows that the solution $t^*$ satisfies
$$ t^* \geq 0.8 \log{(e\lambda/\mu)} \lambda^{-1}.$$
Simultaneously,  $(ex)^{1} + (ex)^{x} < 1$ on $(0, 0.02)$, which shows
$$ t^* \leq  \log{(e\lambda/\mu)} \lambda^{-1} \qquad \mbox{for} \qquad \mu \leq \frac{\lambda}{50}.$$
If $\mu \geq \lambda/50$, then 
$$ e^{- \lambda t} + e^{-\mu t} \leq e^{-\lambda t} + e^{-\lambda t/50} < 1 \qquad \mbox{for some }~t \sim \frac{2.8822}{\lambda} \leq  3\log{ \left(  \frac{e\lambda}{\mu}  \right)} \lambda^{-1}$$
which shows
$$  0.8\log{ \left(  \frac{e\lambda}{\mu}  \right)} \lambda^{-1} \leq t^* \leq  3\log{ \left(  \frac{e\lambda}{\mu}  \right)} \lambda^{-1}.$$
This shows the desired result, the case of Graph Laplacians is completely identical.
\end{proof}

\subsection{Proof of the Corollaries}

\begin{proof}
Let us assume the statement fails and product $\phi_{\mu}\phi_{\lambda}$ 
has a large portion of its $L^2-$energy at frequencies $\leq \lambda/\log{(e\lambda/\mu)}$ 
$$      \sum_{\lambda_k \leq \lambda/\log{(e\lambda/\mu)}}^{}{  |\left\langle \phi_{\mu}\phi_{\lambda}, \phi_{k} \right\rangle|^2} \geq \delta^2 \| \phi_{\mu}\phi_{\lambda} \|^2_{L^2}.$$
Then, clearly,
$$ \| e^{t \Delta} (\phi_{\mu}\phi_{\lambda}) \|_{L^2}^2 = \sum_{k=0}^{\infty}{ e^{-2\lambda_k t}  |\left\langle \phi_{\mu}\phi_{\lambda}, \phi_{k} \right\rangle|^2} \geq 
 \sum_{\lambda_k \leq \lambda/\log{(e\lambda/\mu)}}^{\infty}{ e^{-2\lambda_k t}  |\left\langle \phi_{\mu}\phi_{\lambda}, \phi_{k} \right\rangle|^2}.$$
For $t \sim \log{(e\lambda/\mu)} \lambda^{-1}$, we see that these exponential factors are roughly $\sim 1$ and thus
$$ \| e^{t \Delta} (\phi_{\mu}\phi_{\lambda}) \|_{L^2}^2 \geq e^{-2} \delta^2 \| \phi_{\mu}\phi_{\lambda} \|^2_{L^2},$$
which then implies the desired result since heat evolution and local correlation coincide for that value of $t$.
The proof of Corollary 2 follows the same lines: assume that for $t \sim \log{(e\lambda/\mu)} \lambda^{-1}$
 $$  \left\|  \int_{M}{ p(t,x,y)( \phi_{\lambda}(y) - \phi_{\lambda}(x))( \phi_{\mu}(y) - \phi_{\mu}(x)) dy} \right\|_{L^2_x}  \geq \delta \| \phi_{\mu} \phi_{\lambda} \|_{L^2}.$$
The main result  implies
 $$ \left\| \int_{\Omega}{ p(t,x,y)( \phi_{\lambda}(y) - \phi_{\lambda}(x))( \phi_{\mu}(y) - \phi_{\mu}(x)) dy} \right\|_{L^2_x}^2 =  \| e^{t \Delta} (\phi_{\mu}\phi_{\lambda}) \|_{L^2}^2.$$
At the same time, we can write
\begin{align*}
\delta^2 \| \phi_{\mu}\phi_{\lambda} \|_{L^2}^2  \leq \| e^{t \Delta} (\phi_{\mu}\phi_{\lambda}) \|_{L^2}^2 &=   \sum_{k =1}^{\infty}{   e^{-2\lambda_k t} |\left\langle \phi_{\mu}\phi_{\lambda}, \phi_{k} \right\rangle|^2} \\
&\leq \sum_{\lambda_k \leq c \lambda }^{\infty}{  |\left\langle \phi_{\mu}\phi_{\lambda}, \phi_{k} \right\rangle|^2} +
 \sum_{\lambda_k \geq c \lambda }^{\infty}{ e^{- 2c \lambda_k t} |\left\langle \phi_{\mu}\phi_{\lambda}, \phi_{k} \right\rangle|^2} \\
&\leq \sum_{\lambda_k \leq c \lambda }^{\infty}{ |\left\langle \phi_{\mu}\phi_{\lambda}, \phi_{k} \right\rangle|^2}  + e^{-2 c} \sum_{\lambda_k \geq c \lambda }^{\infty}{  |\left\langle \phi_{\mu}\phi_{\lambda}, \phi_{k} \right\rangle|^2}. 
\end{align*}
This then implies
$$ \sum_{\lambda_k \leq c \lambda }^{\infty}{|\left\langle \phi_{\mu}\phi_{\lambda}, \phi_{k} \right\rangle|^2}  \geq \frac{ \delta^2 e^{2c} - 1}{ e^{2c} - 1}.$$
The estimate is nontrivial for $c \sim \log{(1/\delta)}$ in which case the right-hand side is $\sim \delta^2.$
\end{proof}

\section{ Remarks}
We believe that the local correlation
$$  \int_{M}{ p(t,x,y)( \phi_{\lambda}(y) - \phi_{\lambda}(x))( \phi_{\mu}(y) - \phi_{\mu}(x)) dy} $$
is of intrinsic interest for possibly other reasons as well.  Note that the eigenfunctions
$ \sin{(5 x)}$ and $ \sin{(3x)} \sin{(4y)}$
on $\mathbb{T}^2$ have the same eigenvalue but behave very differently from each other. Conversely, $\sin{(5x)}$
and $\sin{(4 x)}$ correspond to different eigenvalues but are somehow 'more' similar. This question has tremendous importance
in applications and was addressed by Coifman in various collaborations \cite{co1, co2, co3, co4} and was the subject of the
Ph.D. theses of Jerrod Ankenmann \cite{ankenmann}, Sarah Constantin \cite{sarah} and William Leeb \cite{leeb} as well as additional work
of Leeb \cite{leeb2} and Ankenmann \& Leeb \cite{ank2}. A crucial ingredient
is to measure similarity of eigenfunctions by looking at inner products restricted to a subset of the manifold. 
Considering $L^2-$normalized eigenfunctions $e^{in \theta}/\sqrt{2\pi}$, we see that
$$ \left| \int_{a}^{b}{ e^{i(n-n')} d\theta} \right|= \frac{1}{|n-n'|} \left|2-2\cos{(n-n')(b-a)} \right|,$$
which is a quantity that does not depend on $a,b$ but merely on the scale $|b-a|$.  The Ph.D. thesis of Constantin \cite{sarah}
contains a particularly nice example where $\mathbb{S}^3$ is equipped with two different metrics leading to a different 
affinity between the eigenfunctions (that coincide in both).
We note that this idea of local testing is very similar to our global correlation functional.
This raises the question whether the global correlation could be used as an affinity measure in applications.\\

\textbf{Acknowledgments.} The author was asked how $\left\langle \phi_{\mu} \phi_{\lambda}, \phi_{\nu} \right\rangle$ behaves
independently by Martin Ehler during the workshop 'Tractability of High Dimensional Problems and
 Discrepancy' hosted by the Erwin Schr\"odinger Institute (ESI) in Vienna from 25.-29.Sep. 2017 and by Michael Bronstein,
Frederico Monti \& Emanuele Rodola during a visit at Yale in Nov 2017 and their questions inspired this paper. The author is also grateful to Raphy Coifman for
many discussions about the dual geometry induced by eigenfunctions.

\end{document}